\documentclass[12pt,A4paper]{article}
\usepackage{amssymb,amsmath,amsfonts,mathrsfs}

\newtheorem{thm}{Theorem}[section]

\newtheorem{cor}[thm]{Corollary}

\newcommand{\N}{V}

\def\qed{\nopagebreak\hfill{\rule{4pt}{7pt}}\medbreak}

\makeatletter \@addtoreset{equation}{section} \makeatother

\begin{document}
\begin{center}
{\Large\bf On Singletons and Adjacencies of\\[6pt] Set Partitions of Type $B$}
\end{center}

\begin{center}
William Y.C. Chen and David G.L. Wang\\[6pt]
Center for Combinatorics, LPMC-TJKLC\\[6pt]
Nankai University, Tianjin 300071, P.R. China\\[6pt]
chen@nankai.edu.cn, wgl@cfc.nankai.edu.cn
\end{center}

\begin{abstract}
 We show that the joint distribution of the number of singleton pairs and the number
 of adjacency pairs is symmetric over the  set partitions of type $B_n$ without zero-block,
  in analogy with the result of Callan for ordinary partitions.
\end{abstract}



\section{Introduction}\label{Section_Intro}

The main objective of this paper is to give a type $B$ analogue of
 an elegant property of set partitions
discovered by Bernhart \cite{Ber99}, that is,
the number $s_n$ of partitions of $[n]=\{1, 2, \ldots, n\}$ without singletons is equal to the number $a_n$ of partitions of $[n]$ for which no block contains two adjacent elements $i$ and $i+1$ modulo $n$. In fact, it is easy to show that  $s_n$ and $a_n$ have the same formula by the principle of inclusion-exclusion. Bernhart gave a recursive proof of the
fact that $s_n=a_n$ by  showing that $s_n+s_{n+1}=B_n$ and $a_n+a_{n+1}=B_n$, where $B_n$ denotes the Bell number, namely, the number of partitions of $[n]$.
  As noted by Bernhart, there may be no simple way to bring the set of partitions of $[n]$ without singletons and the set of partitions of $[n]$ without adjacencies into a one-to-one correspondence.

 From a different perspective, Callan~\cite{Cal05} found a bijection in terms of an
 algorithm that  interchanges singletons and adjacencies. In fact, Callan has established a stronger statement that the joint distribution of the number of singletons and the number of adjacencies  is symmetric over the set of partitions of $[n]$.  While the proof of Callan is purely combinatorial, we feel that there is still some truth in the remark of Bernhart.

The study of singletons and adjacencies of partitions goes back to Kreweras \cite{Kre72} for
noncrossing partitions. Kreweras has shown that the number of noncrossing partitions of
$[n]$ without singletons equals the number of noncrossing partitions of $[n]$ without
adjacencies. Bernhart~\cite{Ber99} found a combinatorial proof of this assertion. Deutsch and
Shapiro~\cite{DS01} considered noncrossing partitions of $[n]$ without visible singletons
and showed that such partitions are enumerated by the Fine number. Here a visible singleton
of a partition  means a singleton not covered by any arc in the linear
representation. Canfield \cite{Can95} has shown that the average number of singletons in a
partition of $[n]$ is an increasing function of $n$.
Biane~\cite{Bia97} has derived a
bivariate generating function for the number of
partitions of $[n]$ containing a given
number of blocks but no singletons.
Knuth \cite{Knu05} proposed the problem of finding the
generating function for  the number
of partitions of $[n]$ with a given number of blocks
but no adjacencies. The generating function has been
found by several problem solvers.
The sequence of the numbers $s_n$ is listed as
the entry  A000296 in Sloane~\cite{Slo}.

It is natural to wonder  whether there exist a type $B$ analogue
of Bernhart's theorem and a type $B$ analogue of Callan's algorithm.
 We give the peeling and patching algorithm which implies
 the symmetric distribution of the number of singleton pairs and the number of
 adjacency pairs for type $B$ partitions without zero-block.
   Moreover, we can transform the
   bijection to an involution. This
   involution is described in the last section.

\section{The peeling and patching algorithm}

In this section, we give a type $B$ analogue of Callan's symmetric distribution of singletons and adjacencies. Moreover, the algorithm of Callan can  be extended to the type $B$ case.
This type $B$ algorithm will be called the  peeling and patching algorithm.

A partition of type $B_n$ is a partition
$\pi$ of the set $[\pm n]=\{\pm1, \pm2, \ldots, \pm n\}$
such that for any block $B$ of $\pi$, $-B$ is also a
block of $\pi$, and there is at most one block $B$,
called  zero-block,  satisfying $B=-B$,
see Reiner~\cite{Rei97}. We call $\pm i$ a singleton pair
of $\pi$ if $\pi$ contains a block $\{i\}$, and call
$\pm (j,\,j+1)$ an adjacency pair of $\pi$ if $j$ and $j+1$
(modulo $n$) lie in the same block of $\pi$.
Denote the number of singleton pairs (resp. adjacency pairs)
of $\pi$ by $s_{\pi}$ (resp. $a_{\pi}$). For example, let
\begin{equation}\label{eg_pi0}
\pi=\{\pm\{1\},\ \pm\{2\},\ \pm\{3,11,12\},\ \pm\{4,-7,9,10\},\ \pm\{5,6,-8\}\} .
\end{equation}
Then  have $s_{\pi}=2$ and $a_{\pi}=3$.

Denote by $\N_n$ the set of $B_n$-partitions without zero-block.
The following is  the main result of this paper.

\begin{thm}\label{thm_Callan_B}
The joint distribution of the number of singleton pairs and the number of
adjacency pairs is symmetric over $B_n$-partitions without zero-block. In other words, let
\[
P_n(x,y)=\sum_{\pi\in\N_n}x^{s_{\pi}}y^{a_{\pi}},
\]
we have $P_n(x,y)=P_n(y,x)$.
\end{thm}

For example, there are three $B_2$-partitions without zero-block:
\[
\{\pm\{1\},\,\pm\{2\}\},\ \{\pm\{1,2\}\},\ \{\pm\{1,-2\}\}.
\]
So $P_2(x,y)=x^2+y^2+1$. Moreover,
\begin{align*}
P_3(x,y)&=(x^3+y^3)+3xy+3(x+y),\\[5pt]
P_4(x,y)&=(x^4+y^4)+4(x^2y+xy^2)+8(x^2+y^2)+8xy+4(x+y)+7.
\end{align*}

It should be noted that Theorem~\ref{thm_Callan_B} cannot be
deduced from Callan's result for ordinary partitions. The following
consequence is immediate, which is a type $B$ analogue of
Bernhart's observation.

\begin{cor}
The number of $B_n$-partitions without zero-block
and singleton pairs equals the number of
$B_n$-partitions without zero-block and adjacency pairs.
\end{cor}

To prove Theorem \ref{thm_Callan_B}, we shall
provide an algorithm $\psi\colon\N_n\to\N_n$, called the
peeling and patching algorithm,
such that for any $B_n$-partition $\pi$ without zero-block,
$s_{\pi}=a_{\psi(\pi)}$ and $a_{\pi}=s_{\psi(\pi)}$.

In fact, we need a more general setting
to describe the algorithm.
Let $S=\{\pm t_1,\,\pm t_2,\,\ldots,\,\pm t_r\}$ be a
subset of $[\pm n]$, where
$0<t_1<t_2<\cdots<t_r$.
Let $\pi$ be a partition of the set $S$.
We call $\pi$ a {\em symmetric partition} if for any block $B$
of $\pi$, $-B$ is also a block of $\pi$.
We call $\pm t_i$ a {\em singleton pair} of $\pi$ if $\pi$
contains a block $\{t_i\}$, and call $\pm(t_j,t_{j+1})$
an {\em adjacency pair} of $\pi$ if   $t_j$ and $t_{j+1}$
are contained in the same block.
By convention we consider  $t_{r+1}$ as $t_1$.
We call $\pm t_j$ (resp. $\pm t_{j+1}$)
a {\em left-point-pair} (resp. {\em right-point-pair})
if $\pm(t_j,t_{j+1})$ is an adjacency pair.
For the case $r=1$, the partition $\pi=\{\pm\{t_1\}\}$
 contains exactly one singleton pair $\{\pm t_1\}$
 and one adjacency pair $\pm(t_1,t_1)$.

The peeling and patching algorithm $\psi$
consists of the  peeling procedure $\alpha$
and the patching procedure
 $\beta$. During the peeling procedure, at each step
 we take out the singleton pairs and left-point-pairs,
 until there exists neither singleton pairs nor adjacency pairs.
 During  the
  patching procedure, we first interchange the
  roles of singleton pairs and adjacency pairs, then
  put the singleton pairs and
  left-point-pairs back to the partition.
   It should be emphasized that the patching procedure is not
    just  the reverse of the
   peeling procedure.

\noindent{\bf The peeling procedure $\alpha$.}
Given an input partition $\pi$, let $\pi_0=\pi$.
 We extract the set $S_1$ of singleton pairs
 and the set $L_1$ of left-point-pairs (of adjacency pairs)
  from $\pi_0$. Let $\pi_1$ be the remaining partition.
  Now $\pi_1$ is again a type $B$ partition without zero-block.
  So we can extract the set $S_2$ of singleton pairs
  and extract the set $L_2$ of left-point-pairs from $\pi_1$.
Denote by $\pi_2$ be the remaining partition. Repeating this process,
we eventually obtain a partition $\pi_k$ that does not
have any singleton pairs or adjacency pairs. Notice that it is possible
that $\pi_k$ is the empty partition.

For example, consider the partition $\pi$ in (\ref{eg_pi0}), that is,
\[
\pi=\{\pm\{1\},\ \pm\{2\},\ \pm\{3,11,12\},\ %
\pm\{4,-7,9,10\},\ \pm\{5,6,-8\}\}.
\]
 The peeling procedure is illustrated by Table \ref{t1}.

\begin{table}[h]
\centering
\begin{tabular}{|c|c|c|c|}
\hline\rule{0pt}{15pt}$j$
&$S_j$&$L_j$&$\pi_j$\\[3pt]
\hline\rule{0pt}{13pt}1%
&$\pm 1,\, \pm 2$%
&$\pm 5,\, \pm 9,\, \pm 11$%
&$\pm\{3,12\},\ \pm\{4,-7,10\},\ \pm\{6,-8\}$\\[3pt]
\hline\rule{0pt}{15pt} 2%
&$\emptyset$%
&$\pm 12$%
&$\pm\{3\},\ \pm\{4,-7,10\},\ \pm\{6,-8\}$\\[3pt]
\hline\rule{0pt}{15pt} 3%
&$\pm 3$
&$\emptyset$%
&$\pm\{4,-7,10\},\ \pm\{6,-8\}$\\[3pt]
\hline\rule{0pt}{15pt} 4%
&$\emptyset$%
&$\pm 10$%
&$\pm\{4,-7\},\ \pm\{6,-8\}$\\[3pt]
\hline%
\end{tabular}
\caption{\label{t1} The peeling procedure.}
\end{table}

\noindent{\bf The patching procedure $\beta$. }%
Let $\sigma_k=\pi_k$.
As the first step,
we interchange the roles of the singleton-sets $S_i$
and the adjacency-sets represented by $L_i$.
To be precise, we patch the elements of $S_i$ and $L_i$
into the partition $\sigma_i$ which will be obtained
recursively from $\sigma_{i+1}$,
so that $S_i$ (resp. $L_i$) is the right-point-set
(resp. singleton-set) of the resulting partition
$\sigma_{i-1}$.
So $S_i$ represents the set of adjacency pairs of $\sigma_{i-1}$.

We start the patching procedure by putting the elements
of $S_k$ and $L_k$ back to $\sigma_k$ in such a way
that the resulting partition $\sigma_{k-1}$
contains $S_k$ (resp. $L_k$) as its right-point-set
(resp. singleton-set).
The existence of such a partition $\sigma_{k-1}$ will be
confirmed later.
Next, in the same manner we  put the elements
of $S_{k-1}$ and $L_{k-1}$
back into $\sigma_{k-1}$ to get $\sigma_{k-2}$.
Repeating this process, we finally arrive at a
partition $\sigma_0$, which is defined to be the output of the
patching procedure.

Now let us describe the process of
constructing $\sigma_{k-1}$.
Suppose that the underlying set of $\pi_{k-1}$
is $\{\pm t_1,\,\pm t_2,\,\ldots,\,\pm t_r\}$,
where $0<t_1<t_2<\cdots<t_r$.

Consider the case that  $\sigma_k(=\pi_k)$ is the empty partition.
The last step of the peeling procedure implies that
$\pi_{k-1}$ must be of special form, namely,
either there is only one block in $\pi_{k-1}$,
or every block of $\pi_{k-1}$ contains exactly one element.
Define $\sigma_{k-1}$ to be
$\{\pm\{t_1\},\ \pm\{t_2\},\ \ldots,\ \pm\{t_r\}\}$
if there is only one block in $\pi_{k-1}$; otherwise, set
$\sigma_{k-1}=\{\pm\{t_1,\,t_2,\,\ldots,\,t_r\}\}$.
When  $r=1$, it is clear to see that $\sigma_{k-1}$ is well-defined.

We now assume that $\sigma_k$ is not empty.
We can uniquely decompose the set $S_t$ into maximal
consecutive subsets of the form
\begin{equation}\label{eq2}
\{\pm t_{i+1},\,\pm t_{i+2},\,\ldots,\,\pm t_{i+h}\}.
\end{equation}
The number of such subsets is at least two.
By the maximality, the element $t_i$ does not appear in $S_k$.
On the other hand, it is clear that $t_i\not\in L_k$ by
the definition of $L_k$. Thus $t_i$ is contained in $\sigma_k$.
This observation allows us to put the elements
$t_{i+1},\,t_{i+2},\,\ldots,\,t_{i+h}$
into the block of $\sigma_k$ containing $t_i$.
Accordingly, we  put $-t_{i+1},\,-t_{i+2},\,\ldots,\,-t_{i+h}$
into the block containing $-t_i$.
After having processed all maximal consecutive subsets of $S_k$,
we put  each element in $L_k$
as a singleton block into the partition $\sigma_k$.
The resulting partition is defined to be $\sigma_{k-1}$.

This completes the description of the step of
constructing $\sigma_{k-1}$. Since $\sigma_k(=\pi_k)$
contains neither singleton pairs nor adjacency pairs,
it is easy to check that $L_k$ (resp. $S_k$)
is the set of singleton pairs (right-point-pairs) of $\sigma_{k-1}$.

For example,  Table \ref{t2} is an illustration of the
patching procedure for partition generated
 in Table \ref{t1}.
\begin{table}[h]
\centering
\begin{tabular}{|c|c|c|c|}
\hline\rule{0pt}{15pt}%
$j$&$S_j$&$L_j$&$\sigma_j$\\[3pt]
\hline\rule{0pt}{15pt} 4%
&$\emptyset$%
&$\pm 10$%
&$\pm\{4,-7\},\ \pm\{6,-8\}$\\[3pt]
\hline\rule{0pt}{15pt} 3%
&$\pm 3$%
&$\emptyset$%
&$\pm\{4,-7\},\ \pm\{6,-8\},\ \pm\{10\}$\\[3pt]
\hline\rule{0pt}{15pt} 2%
&$\emptyset$%
&$\pm 12$%
&$\pm\{4,-7\},\ \pm\{6,-8\},\ \pm\{3,10\}$\\[3pt]
\hline\rule{0pt}{15pt} 1%
&$\pm 1,\, \pm 2$%
&$\pm 5,\, \pm 9,\, \pm 11$%
&$\pm\{4,-7\},\ \pm\{6,-8\},\ \pm\{3,10\}, \pm\{12\}$\\[3pt]
\hline%
\end{tabular}
\caption{\label{t2} The patching procedure.}
\end{table}
In the last step, patching $S_1$ and $L_1$ to $\sigma_1$,
we finally obtain
\begin{equation}\label{sigma0}
\sigma_0=\{\pm\{1,2,12\},\ \pm\{3,10\},\ \pm\{4,-7\},\ \pm\{5\},\ %
\pm\{6,-8\},\ \pm\{9\},\ \pm\{11\}\}.
\end{equation}

The peeling and patching algorithm $\psi$ is defined by
\[  \psi(\pi)=\beta(\alpha(\pi)) \]
  for any $B_n$-partition $\pi$ without zero-block.
Keep in mind that there is a step of
interchanging the roles  of singleton pairs
and adjacency pairs at the beginning of the patching procedure.
We are now ready to give a proof of Theorem \ref{thm_Callan_B}.

\noindent{\it Proof of Theorem \ref{thm_Callan_B}. }
We aim to show that the peeling and patching algorithm $\psi$
gives a bijection on $B_n$-partitions without zero-block,
 which interchanges the number of singleton pairs and the number
 of adjacency pairs.

It is easy to see that the inverse algorithm
can be described as follows. It is in fact the composition of
another peeling procedure and another patching procedure.
To be precise, let $\sigma$ be the input partition.
Let $\sigma_0=\sigma$. We first peel the singleton pairs and
right-point-pairs at each step,
until we obtain a partition $\sigma_k$
which has neither singleton pairs nor adjacency pairs. Then,
based on the partition $\pi_k=\sigma_k$,
we recursively patch the elements
that have been taken out before.
Meanwhile,  we also need to interchange
the roles of
the singleton-sets and right-point-sets
at the beginning of this patching procedure.
Finally,  we get a
partition, as the output of the inverse algorithm.
Therefore, $\psi$ is a bijection which exchanges the
number of singleton pairs and the number of adjacency pairs.
 This completes the proof. \qed

An illustration of the peeling and patching algorithm is
given by (\ref{eg_pi0}),
Table \ref{t1}, Table \ref{t2},
and (\ref{sigma0}).

To conclude this section,
we give the generating function for the
number $s_n^B$ of $B_n$-partitions without zero-block and singleton
pairs, that is,
\begin{equation}\label{EGF}
\sum_{n\ge0}s_n^B\frac{x^n}{n!}
=\exp\left(\sinh(x)e^x-x\right).
\end{equation}
By the principle of inclusion-exclusion, we obtain
\begin{equation}\label{eq1}
s_n^B=\sum_{k=0}^n(-1)^{n-k}{n\choose k}\sum_{j=0}^k2^{k-j}S(k,j),
\end{equation}
where $S(k,j)$ is the Stirling number of the second kind, and
$2^{k-j}S(k,j)$ is the number of partitions in $\N_k$
containing exactly $2j$ blocks.
The formula~(\ref{EGF}) can be easily
derived from~(\ref{eq1}).

\section{From bijection to involution}

The bijection given in the previous section is not an involution although it
interchanges the number of  singleton pairs and the number of adjacency pairs.
In this section,  we show that the peeling and patching algorithm can be turned
into an involution. Such an involution  for ordinary partitions
 has been given by Callan \cite{Cal05}.

For any $i\in [n]$, we define the {\em complement} of $i$ to be $n+1-i$, and the complement of $-i$ as  $-(n+1-i)$. This notion can be extended naturally to any symmetric partition $\pi$ of $[\pm n]$ by taking the complement for each element in the partition. The complement of $\pi$
is denoted by $\omega(\pi)$. It is clear that  $\omega$ is an involution.
Assume that $\sigma_0$ is given  in~(\ref{sigma0}). We have
\begin{equation}\label{omega}
\omega(\sigma_0)=\{\pm\{1,11,12\},\ \pm\{2\},\ \pm\{3,10\},\ \pm\{4\},\
\pm\{5,-7\},\ \pm\{6,-9\},\ \pm\{8\}\}.
\end{equation}

In light of the complementation operation,
we get an involution based on the
peeling and patching algorithm. The proof is a straightforward verification and hence is omitted.

\begin{thm}
The mapping  $\omega\circ\psi$ is an involution on $B_n$-partitions without zero-block, which
interchanges the number of singleton pairs and the number of adjacency pairs.
\end{thm}

Let us give an example to demonstrate that $\omega\circ \psi$ is involution, that is
\begin{equation}\label{eq5}
\omega(\psi(\pi))=\psi^{-1}(\omega(\pi)).
 \end{equation}
 Consider the partition $\pi$ in (\ref{eg_pi0}). In this case, the left hand side of (\ref{eq5}) is $\omega(\sigma_0)$ in (\ref{omega}). On the other hand,
\[
\omega(\pi)=\{\pm\{1,\,2,\,10\},\ \pm\{3,\,4,\, -6,\,9\},\ \pm\{5,\,
-7,\,-8\},\ \pm\{11\},\ \pm\{12\}\}.
\]
Applying the procedure $\beta^{-1}$, we obtain the Table \ref{t3},
where $R_j$ (resp. $S_j$) denotes the set of right-point-pairs
(singleton pairs). Next, by the procedure $\alpha^{-1}$,
we get the Table \ref{t4}.
Finally, putting  $R_1$ and $S_1$ back to $\pi_1$,
we arrive at the partition $\pi_0$
which is in agreement with (\ref{omega}).

\begin{table}[h]
\centering
\begin{tabular}{|c|c|c|c|}
\hline\rule{0pt}{15pt}
$j$&$R_j$&$S_j$&$\sigma_j$\\[3pt]
\hline\rule{0pt}{15pt} 1%
&$\pm 2,\, \pm 4,\, \pm 8$%
&$\pm 11,\, \pm 12$%
&$\pm\{1,10\},\ \pm\{3,-6,9\},\ \pm\{5,-7\}$\\[3pt]
\hline\rule{0pt}{15pt} 2%
&$\pm 1$%
&$\emptyset$%
&$\pm\{10\},\ \pm\{3,-6,9\},\ \pm\{5,-7\}$\\[3pt]
\hline\rule{0pt}{15pt} 3%
&$\emptyset$%
&$\pm 10$%
&$\pm\{3,-6,9\},\ \pm\{5,-7\}$\\[3pt]
\hline\rule{0pt}{15pt} 4%
&$\pm 3$%
&$\emptyset$%
&$\pm\{5,-7\},\ \pm\{6,-9\}$\\[3pt]
\hline%
\end{tabular}
\caption{\label{t3} The procedure $\beta^{-1}$.}
\end{table}

\begin{table}[h]
\centering
\begin{tabular}{|c|c|c|c|}
\hline\rule{0pt}{15pt}%
$j$&$R_j$&$S_j$&$\pi_j$\\[3pt]
\hline\rule{0pt}{15pt} 4%
&$\pm 3$%
&$\emptyset$%
&$\pm\{5,-7\},\ \pm\{6,-9\}$\\[3pt]
\hline\rule{0pt}{15pt} 3%
&$\emptyset$%
&$\pm 10$%
&$\pm\{3\},\ \pm\{5,-7\},\ \pm\{6,-9\}$\\[3pt]
\hline\rule{0pt}{15pt} 2%
&$\pm 1$%
&$\emptyset$%
&$\pm\{3,10\},\ \pm\{5,-7\},\ \pm\{6,-9\}$\\[3pt]
\hline\rule{0pt}{15pt} 1%
&$\pm 2,\,\pm 4,\,\pm 8$%
&$\pm 11,\, \pm 12$%
&$\pm\{1\},\ \pm\{3,10\},\ \pm\{5,-7\},\ \pm\{6,-9\}$\\[3pt]
\hline%
\end{tabular}
\caption{\label{t4} The procedure $\alpha^{-1}$.}
\end{table}

\noindent{\bf Acknowledgments.} This work was supported by the 973
Project, the PCSIRT Project of the Ministry of Education, and the
National Science Foundation of China.

\end{document}